\documentclass[a4paper,10pt]{article}
\usepackage{a4,ae,amsmath,amssymb,verbatim}   

\newtheorem{Satz}{Proposition} 
\newtheorem{Lemma}[Satz]{Lemma} 
\newtheorem{Rem}[Satz]{Remark} 
\newtheorem{Cor}[Satz]{Corollary}
\newtheorem{Def}[Satz]{Definition} 
\newtheorem{Thm}[Satz]{Theorem} 
\newtheorem{BSP}{Example}
\newenvironment{Bsp}{\begin{BSP}\normalfont}{\end{BSP}}

\newcommand{\PROOF}[1][]{\textsc{Proof#1:} }
\newcommand{\END}{\hfill$\square$\smallskip} 
\newcommand{\ENDCLAIM}{\hfill$\lozenge$\smallskip} 

\newcommand{\Emb}{\mathrm{Emb}}
\newcommand{\End}{\mathrm{End}}
\newcommand{\Sym}{\mathrm{Sym}}
\newcommand{\rck}[1]{\underline{#1}}
\newcommand{\mc}{\mathcal}
\newcommand{\mf}{\mathfrak}
\newcommand{\Aut}{\mathrm{Aut}}
\newcommand{\verum}{\top}
\newcommand{\falsum}{\bot}
\newcommand{\id}{\mathrm{id}}
\newcommand{\conc}{\text{\hspace*{-0.25em}\t{}\ }}                              
\newcommand{\surj}{\twoheadrightarrow}                  
\newcommand{\isom}{\cong}                       
\newcommand{\ohne}{\setminus}

\newenvironment{LList}[1][\alph]{\begin{list}{{\rm\bf(#1{lemmaitem1})}}{\usecounter{lemmaitem1}\itemsep1ex\setlength{\topsep}{0.5\itemsep}\parsep0ex\labelsep1ex\settowidth{\labelwidth}{(m)}\setlength{\leftmargin}{\labelwidth}\addtolength{\leftmargin}{\labelsep}}}{\end{list}}

\begin{document} 

\title{$\aleph_0$-categorical Structures: Endomorphisms and Interpretations }
\author{Manuel Bodirsky%
   \footnote{Laboratoire d'Informatique de l'\'Ecole Polytechnique,
     91128 Palaiseau, France,
     bodirsky@lix.polytechnique.fr}
   \addtocounter{footnote}{5},
   Markus Junker%
   \footnote{Albert--Ludwigs--Universit{\"a}t Freiburg,
     Abt. f{\"u}r Math. Logik, Eckerstra{\ss}e 1, 79102
     Freiburg, Germany, markus.junker@math.uni-freiburg.de}}
\maketitle

\begin{abstract}
\noindent
  We extend the Ahlbrandt--Ziegler analysis of interpretability in
  $\aleph_0$-categorical structures by showing that existential
  interpretation is controlled by the monoid of self--embeddings and
  positive existential interpretation of structures without constant
  endomorphisms is controlled by the monoid of endomorphisms in the
  same way as general interpretability is controlled by the
  automorphism group. 
\end{abstract}

\section{Introduction}

$\aleph_0$-categorical structures (often called $\omega$-categorical
structures) appear quite naturally in mathematics, and have
extensively been studied by model theorists. 
They appear for example as countable universal structures for classes
of finite structures with the amalgamation property. The best known
example might be the countable random graph, which can be seen as a
universal amalgam of the class of all finite graphs. The 
$\aleph_0$-categorical structures can also be
characterised by a transitivity property of their automorphisms groups,
which are so-called ``oligomorphic permutation groups'', and therefore
they are also interesting for and have been studied by group theorists. 
More on $\aleph_0$-categorical structures can for example be found in
\cite{Hodges}, Sections 7.3 and 7.4, \cite{KM} and \cite{Cameron}.

In fact, much of an $\aleph_0$-categorical structure is coded in
its automorphism group. Ahlbrandt and Ziegler in \cite{AZ} have shown that 
a countable $\aleph_0$-categorical structure is, up to bi-interpretability, 
determined by its automorphism group as a topological group. We extend this 
analysis and show that, with certain unavoidable restrictions, existential 
interpretability is controlled by the monoid of self-embeddings and 
positive existential interpretability by the endomorphism monoid. 

It would be interesting to further extend the theory (as far as
possible) to primitive positive interpretability on the one hand and
polymorphism clones on the other hand; a characterisation of
primitive positive interpretability in terms of the topological
polymorphism clone would have interesting consequences for the study
of the computational complexity of constraint satisfaction problems in
theoretical computer science.

\section{Endomorphisms}
\subsection{Preservation theorems for $\aleph_0$-categorical theories}

In this paper, we only consider structures $\mf M$ in a countable
signature without function symbols (i.e. relational possibly with
constants). We denote by $\Aut(\mf M)$ the automorphism group of $\mf
M$, by $\Emb(\mf M)$ the monoid of embeddings%
\footnote{i.e. isomorphisms onto a substructure, or equivalently
  strong injective homomorphisms}  
of $\mf M$ into $\mf M$, and by $\End(\mf M)$ the monoid of all
endomorphisms of $\mf M$. Then $\Aut(\mf M) \subseteq \Emb(\mf M)
\subseteq \End(\mf M) \subseteq \mbox{${}^M\!M$}$.
All these monoids carry the topology of pointwise convergence, a basis
of open neighbourhoods of which is given by the sets
$U_{\bar a,\bar b} = \{\sigma \mid \bar a^\sigma = \bar b\}$. 
Finally, $\Sym(M)$ denotes the symmetric group on $M$.

\begin{Rem} 
\textbf{\upshape{(a)}}
$\Emb(\mf M)$ and $\End(\mf M)$ are closed in \mbox{${}^M\!M$},
because if a map is not a homomorphisms, not injective or not strong,
then this is already witnessed by a finite tuple, hence a complete open
neighbourhood does lack this property. $\Aut(\mf M)$ is closed in 
$\Sym(M)$, and more generally in the set
of all surjections $M \to M$, but in general not in \mbox{${}^M\!M$}.  

\textbf{\upshape{(b)}}
$\Aut(\mf M) = \Emb(\mf M) \cap \Sym(M)$, but in general there are
more bijective endomorphisms than automorphisms as they need not to be
strong. But then their inverse maps are not homomorphisms. It follows
that $\Aut(\mf M)$ equals the set of invertible elements of $\End(\mf M)$ 
and therefore the largest subgroup of $\End(\mf M)$.
\end{Rem}

\textbf{Notations and conventions:} 
For the sake of this paper, we call a structure \emph{$\aleph_0$-categorical} 
if it is finite or a countable model of an $\aleph_0$-categorical theory with
an at most countable language. We will freely use the characterisation
of Engeler, Ryll--Nardzewski and Svenonius (see \cite{Hodges},
Theorem~7.3.1), which in particular implies the \emph{ultrahomogeneity} of
an $\aleph_0$-categorical structure: any two tuples of same type are
conjugate under the automorphism group. 

We let endomorphisms act from the right side and write $x^\sigma$ for
$\sigma(x)$ and $x^{\sigma\tau}$ for $\tau(\sigma(x))$, and in
particular $\sigma\tau$ for $\tau\circ\sigma$. 

Formulae, definability etc, are meant without parameters, unless
otherwise specified. For the present paper, it doesn't make a difference
whether we understand ``existential formula'' and ``positive formula''
up to logical equivalence or not. It is a classical result that this
works as well for ``positive existential'', i.e. positivity and
existentiality can be realised simultaneously (see e.g.\ \cite{CK}
Exercise 5.2.6). 

Let $\Sigma$ be a set of maps from $M$ to $M$. Each $\sigma \in \Sigma$ induces a 
map $M^k \to M^k$ (acting component by component), also denoted by $\sigma$. 
The \emph{orbit} of $\bar m \in M^k$ under $\Sigma$ is the set of images 
$\bar m^\Sigma = \{\bar m^\sigma \mid \sigma \in \Sigma\}$. In general, the
orbits are not the classes of an equivalence relation. If $X \subseteq M^k$, 
then $X$ is called \emph{closed under $\Sigma$} if for all $x \in X$, the 
orbit $x^\Sigma$ is contained in $X$.

\begin{Satz}[(a),(b) in \cite{BP}, Theorem~5] \label{S-BP}
Let $M$ be an $\aleph_0$-categorical structure, and $X \subseteq M^k$. 
\begin{LList}
\item 
  $X$ is existentially definable in $\mf M$ if and only if $X$
  is closed under $\Emb(\mf M)$. 
\item 
  $X$ is positive existentially definable in $\mf M$ if and only if $X$
  is closed under $\End(\mf M)$.  
\item 
  $X$ is positively\/%
 \footnote{$\verum$ and $\falsum$ are positive formulae.} 
   definable in $\mf M$ if and only if $X$ is closed under all
   surjective endomorphisms of $\mf M$. 
\item 
  $X$ is positive existentially definable in $\mf M$ in the language with 
  $\neq$ if and only if $X$ is closed under all injective
  endomorphisms of $\mf M$. 
\item  
  $X$ is positively definable in $\mf M$ in the language with 
  $\neq$ if and only if $X$ is closed under all bijective
  endomorphisms of $\mf M$. 
\end{LList}
\end{Satz} 

\PROOF
It is well known by Ryll--Nardzewski etc. that $X$ is definable if and
only if $X$ is invariant under $\Aut(\mf M)$, and because $\Aut(\mf M)$ 
is a group, this is equivalent to  being closed under $\Aut(\mf M)$. 
Therefore, we may assume that $X$ is definable by a formula $\phi$.

(b) and (c):
If $\mf M \models \phi(\bar a)$ and $\phi$ is positive, then 
$\sigma(\mf M)\models \phi(\bar a^\sigma)$ for every homomorphism $\sigma$. 
If $\sigma$ is surjective or if $\phi$ is in addition existential
modulo $T$, then it follows that $\mf M\models\phi(\bar a^\sigma)$.  

For the other direction, we need the classical {\L}os--Tarski and Lyndon 
preservation theorems (see \cite{Hodges} Theorem 6.5.4 and Corollary 10.3.5). 
By these well-known theorems, if $\phi$ is not positive existential
(not positive), then there are models $\mf M_i \models T$, a
(surjective) homomorphism $\sigma: \mf M_1 \to \mf M_2$ and $\bar a$
in $\mf M_1$ with $\mf M_1 \models \phi(\bar a)$ and 
$\mf M_2 \not\models \phi(\bar a^{\sigma})$. Now choose a countable elementary 
substructure of $(\mf M_1,\mf M_2,\sigma,\bar a)$. Up to isomorphism, it has 
the form $(\mf M,\mf M,\sigma',\bar a)$, where $\sigma'$ is a (surjective)
endomorphism of $\mf M$. Then we get $\mf M\models \phi(\bar a)$, but 
$\mf M \not\models \phi(\bar a^{\sigma'})$. 

(d) follows from (b) and (e) from (c) just by adding $\neq$ to the language. 
In the same way (a) follows from (d) by adding negations of all the basic 
relations to the language. 
\END

Clearly, a set $X\subseteq M^k$ is closed under $\Sigma \subseteq {}^M\!M$ if
and only if $X$ is closed under the closure of $\Sigma$ in ${}^M\!M$. 
Therefore, if $\Sigma_1$ and $\Sigma_2$ are dense in each other, then 
syntactical properties characterised by $\Sigma_1,\Sigma_2$ are equivalent. 
If $\Sigma_1 = \Aut(\mf M) \subseteq \Sigma_2$, then the converse also holds, 
which we will prove for the example $\Sigma_2 = \Emb(\mf M)$:

\begin{Cor}[Bodirsky, Pinsker] 
  An $\aleph_0$-categorical theory with countable model $\mf M$ is model 
  complete if and only if $\Aut(\mf M)$ is dense in $\Emb(\mf M)$. 
\end{Cor} 

\PROOF 
``$\Longleftarrow$'' follows from the general remarks above. 

``$\Longrightarrow$'': According to Proposition~\ref{S-BP}\,(a), 
self-embeddings preserve existential types, hence complete types in case the  
theory is model complete. This implies $\bar a^\sigma \equiv\bar a$ for all 
finite tuples $\bar a$ in $\mf M$ and all $\sigma \in \Emb(\mf M)$. By
the ultrahomogeneity of an $\aleph_0$-categorical model (tuples of
same type are conjugate under the automorphism group), this is
equivalent to $\Aut(\mf M)$ being dense in $\Emb(\mf M)$.
\END

In the same style, every definable set is positively definable in $\mf M$, if 
all surjective homomorphisms are automorphisms.

\subsection{Topology}

Let $\mf M$ be an $\aleph_0$-categorical structure and $T$ its theory. 
We consider the topological space $\End(\mf M)/\Aut(\mf M)$ of \emph{right 
cosets} of $\Aut(\mf M)$ in $\End(\mf M)$, i.e. the quotient of $\End(\mf M)$ 
by the equivalence relation 
$$\sigma \sim \sigma' \iff \text{ there is } \alpha\in\Aut(\mf A) 
  \text{ with } \sigma' = \sigma\alpha,$$ 
equipped with the quotient topology, the finest topology which turns 
$\pi: \sigma\mapsto\sigma/_\sim$ into a continuous map. Inverse images of the 
open sets are open sets in $\End(\mf M)$ of the form $X \cdot\Aut(\mf M)$
for open $X \subseteq \End(\mf M)$.

\begin{Lemma} \label{L-kompakt}
  $\End(\mf M)/\Aut(\mf M)$ and $\Emb(\mf M)/\Aut(\mf M)$ are compact. 
\end{Lemma} 

\PROOF 
As $\Emb(\mf M)$ is closed in $\End(\mf M)$ and a union of right cosets of 
$Aut(\mf M)$, it is sufficient to show the first claim. Consider an open 
covering $(U_i)_{i\in I}$ of $\End(\mf M)/\Aut(\mf M)$. We may assume that the 
inverse images $\tilde U_i := \pi^{-1}[U_i]$ in $\End(\mf M)$ are of the form 
$U_{\bar c_i,\bar d_i}\cdot\Aut(\mf M) = \{\sigma\in \End(\mf M) \mid 
\bar c_i^{\,\sigma} \equiv \bar d_i\}$. Thus the $\tilde U_i$ from an open 
covering of $\End(\mf M)$ by sets which are unions of right cosets. It is 
sufficient to show that $\End(\mf M)$ is covered by finitely many of the 
$\tilde U_i$. Fix an enumeration $(m_i)_{i \in \omega}$ of $M$. If $p$ is an 
$n$-type of $T$, let $U_{p} := U_{(m_0,\dots,m_{n-1}),\bar a}\cdot\Aut(\mf M)$
where $\bar a$ is some/any realisation of $p$. Note that if 
$U_{p}\neq\emptyset$ and $\bar a\models p$, then $(m_0,\dots,m_{n-1}) \mapsto 
\bar a$ is a partial endomorphism. Finally, let us say that an open set $O$ is 
``covered'' if there is an $i\in I$ with $O \subseteq \tilde U_i$. 

If $U_{p}$ is covered for some $n\in\omega$ and each of the finitely many 
$n$-types $p$, then the coverings set form an open sub-covering of 
$\End(\mf M)$. Therefore, we may assume that for each $n$, there is an 
$n$-type $p_n(x_0,\dots,x_{n-1})$ such that $U_{p_n}$ is not covered 
(and in particular, $U_{p_n}\neq\emptyset$). The types $p_n$ form an infinite 
tree under inclusion, which is finitely branched because of the 
$\aleph_0$-categoricity. Hence, by K{\"o}nig's Lemma, there is an infinite 
branch $(p_n)_{n \in \omega}$. If $(a_n)_{n \in \omega}$ realises 
$\bigcup_{n \in \omega}p_n$, then $\sigma: m_n \mapsto a_n$ defines an 
endomorphism of $\mf M$.

Now choose $i$ such that $\sigma \in \tilde U_i$, and let $n$ be big enough 
such that $\bar c_i$ is contained in $\bar m := (m_0,\dots,m_{n-1})$. Then 
$U_{p_n} = U_{\bar m,\bar m^\sigma} \subseteq \tilde U_i$:
contradiction. This shows quasi-compactness.

If $\sigma \not\sim \sigma'$, then there is a tuple $\bar a$ with 
$\bar a^\sigma \not\equiv \bar a^{\sigma'}$. Thus the open neighbourhoods 
$U_{\bar a,\bar a^\sigma}\cdot \Aut(\mf M)$ and 
$U_{\bar a,\bar a^{\sigma'}}\cdot \Aut(\mf M)$ separate $\sigma$ and $\sigma'$.
\END

\begin{Rem} 
  $\Aut(\mf M)$ is not ``normal'' in $\Emb(\mf M)$, i.e. the left coset 
  $\sigma\cdot\Aut(\mf M)$ is in general different from the right coset 
  $\Aut(\mf M)\cdot\sigma$. 
\end{Rem}

\begin{Bsp}
Let $\mf M$ be an equivalence relation with two classes, both countably 
infinite; $\alpha$ is an automorphism that exchanges both classes, and 
$\sigma$ is an embedding that is the identity on one class and non surjective 
on the other class. Then $\sigma^{-1}\alpha\sigma$ can't be extended to an 
automorphism of $\mf M$, i.e. $\alpha\sigma$ is not of the form 
$\sigma\alpha'$ for some automorphism $\alpha'$. 
\end{Bsp}

\section{Interpretations}

The classical theory of interpretations of $\aleph_0$-categorical theories as 
developed by Ahlbrandt and Ziegler in \cite{AZ} is briefly as follows. 
(An account of the theory and more about interpretations can be found
in Section~1 of \cite{KM} and in Section~5 of \cite{Hodges}). 
In \cite{AZ}, $\aleph_0$-categorical structures are considered as a
category with interpretations as morphisms, and ``$\Aut$'' is made
into a functor into the category of topological groups with continuous 
group homomorphisms, where $\Aut(i)$ for an interpretation $i$ of $\mf B$ 
in $\mf A$ is the natural map $\Aut(\mf A) \to \Aut(\mf B)$ induced by $i$. 

\medskip
\textbf{Theorem~1.2 in \cite{AZ}}
\textit{A continuous group homomorphism $f: \Aut(\mf A) \to \Aut(\mf B)$ is of 
  the form $\Aut(i)$ for an interpretation $i$ of $\mf B$ in $\mf A$ if and 
  only if $B$ is covered by finitely many orbits under the image of $f$.}

\medskip
Two interpretations $i_1,i_2$ of $\mf B$ in $\mf A$ are called homotopic if 
$\{(\bar x,\bar y) \mid i_1(\bar x) = i_2(\bar y)\}$ is definable in $\mf A$. 
Two structures $\mf A,\mf B$ are bi-interpretable if there are mutual 
interpretations $i$ and $j$ such that $i \circ j$ and $j \circ i$ are 
homotopic to the identity interpretations $\id_{\mf A}, \id_{\mf B}$ 
respectively.

\medskip
\textbf{Theorem~1.3 and Corollary~1.4 in \cite{AZ}}
\textit{Two interpretations $i_1,i_2$ of $\mf B$ in $\mf A$ are homotopic if 
  and only if $\Aut(i_1) = \Aut(i_2)$. The structures are bi-interpretable if 
  and only if there automorphism groups are isomorphic as topological groups.}

\begin{Rem}
In Theorem~1.2 of \cite{AZ}, one could as well have considered a continuous 
monoid homomorphism $\Aut(\mf A) \to \End(\mf B)$ instead of a continuous 
group homomorphism $\Aut(\mf A) \to \Aut(\mf B)$. This is because a monoid
homomorphism defined on a group is a group homomorphism, and thus 
the group $\Aut(\mf A)$ has to be mapped into the largest group contained in 
$\End(\mf B)$ which is $\Aut(\mf B)$. 
\end{Rem}

Our aim is to extend the classical results to endomorphisms on the one hand 
and to syntactically restricted interpretations on the other hand.

\subsection{The existential case}

Let us call \emph{basic sets of a structure} the universe, the diagonal, the 
interpretations of the relational symbols in the language and the graphs of 
the interpretations of the functions symbols in the language. An 
interpretation of a structure $\mf N$ in a structure $\mf M$ is existential 
(positive existential) if all inverse images of basic sets of $\mf N$ are 
existentially (positive existentially) definable in $\mf M$.

\begin{Thm} \label{Th-exist} 
  Let $\mf A$ be an $\aleph_0$-categorical structure with at least two
  elements. Then $\mf B$ is 
  existentially interpretable in $\mf A$ if and only if there is a continuous 
  monoid homomorphism $f: \Emb(\mf A) \to \End(\mf B)$ such that $B$ is 
  covered by finitely many orbits under the image of $f$, or, equivalently, 
  such that $B$ is covered by finitely many orbits under $f[\Aut(\mf A)]$.
\end{Thm}

\PROOF 
If $B$ is covered by finitely many orbits under $f[\Aut(\mf A)]$, then it is 
also covered by finitely many orbits under $f[\Emb(\mf A)]$. We are going to 
show ``$\Longleftarrow$'' with the weaker and ``$\Longrightarrow$'' with the 
stronger of the two covering conditions. 

``$\Longleftarrow$'': Choose $\bar b = (b_1,\dots,b_k)$ with $b_i \in B$ such
that $B$ is covered by the orbits of the $b_i$ under $f[\Emb(\mf A)]$. 

\textsc{Claim:} There is a finite tuple $\bar a$ in $A$ with the following 
property: If $\bar a^\sigma = \bar a^\tau$ for $\sigma,\tau \in \Emb(\mf A)$, 
then $\bar b^{f(\sigma)} = \bar b^{f(\tau)}$. 

\PROOF[ of the Claim]
We call a tuple $\bar a$ \emph{good for $\sigma$} if $\bar a^\sigma = 
\bar a^\tau$ implies $\bar b^{f(\sigma)} = \bar b^{f(\tau)}$ for all $\tau$. 
Fix $\sigma_0 \in \Emb(\mf A)$. Because $f$ is continuous, $f^{-1}[U_{\bar b, 
\bar b^{f(\sigma_0)}}]$ is an open set containing $\sigma_0$ and thus contains 
a basic open neighbourhood $U_{\bar c,\bar c^{\sigma_0}}$ of $\sigma_0$. Then
$\bar c$ is good for $\sigma_0$ because if $\bar c^{\sigma_0}=\bar c^\tau$, 
then $\tau \in U_{\bar c,\bar c^{\sigma_0}}$, hence $f(\tau) \in U_{\bar b,
\bar b^{f(\sigma_0)}}$ and thus $\bar b^{f(\sigma_0)} = \bar b^{f(\tau)}$. 

Note that $\bar c$ clearly is good for each other $\sigma\in U_{\bar c,
\bar c^{\sigma_0}}$, and also for all $\sigma_0\alpha$ with $\alpha \in 
\Aut(\mf A)$, i.e. for the whole neighbourhood $U_{\bar c,\bar c^{\sigma_0}} 
\cdot \Aut(\mf A)$. For suppose $\bar c^{\sigma_0\alpha} = \bar c^\tau$, then
$\bar c^{\sigma_0} = \bar c^{\tau\alpha^{-1}}$, hence $\bar b^{f(\sigma_0)} = 
\bar b^{f(\tau\alpha^{-1})} = \bar b^{f(\tau)f(\alpha)^{-1}}$  because $f$ is 
a monoid homomorphism and thus maps automorphisms onto automorphisms.
Finally $\bar b^{f(\sigma_0\alpha)} = \bar b^{f(\tau)}$ follows. 

Now we have found a $\bar c_i$ for each $\sigma_i$ which is good for the 
neighbourhood $U_i := U_{\bar c_i^{},\bar c_i^{\sigma_i}} \cdot \Aut(\mf A)$. 
By the compactness of $\Emb(\mf A)/\Aut(\mf A)$ shown in Lemma~\ref{L-kompakt},
finitely many of these neighbourhoods, say $U_1,\dots,U_l$, cover 
$\Emb(\mf A)$. Then  $\bar a:= \bar c_1\conc\cdots\conc \bar c_l$ is a tuple 
which is good  for all $\Emb(\mf A)$. 
\ENDCLAIM

Let $a_1',\dots,a_k'$ be arbitrary pairwise distinct elements of $A$
(or, if $|A| < k$, of some sufficiently large power of $A$). We may
assume that the $a_i'$ appear in the tuple $\bar a$ (otherwise extend
$\bar a$ by the $a_i'$). Now we can continue as in Ahlbrandt--Ziegler: 

\medskip
\textsc{Definition:} \\
Let $U := \big\{(a'_i,\bar a)^\sigma \mid i=1,\dots,k,\, \sigma \in 
\Emb(\mf A)\big\}$ and define $\rck{f}: U \to B$ by $(a'_i,\bar a)^\sigma
\mapsto b_i^{f(\sigma)}$. 

Note that by definition, $U$ is closed under $\Emb(\mf A)$, hence 
existentially definable after Proposition~\ref{S-BP}.

\medskip
\textsc{Claim:}
$\rck{f}$ is well defined and surjective.

If $(a'_i,\bar a)^\sigma = (a'_j,\bar a)^\tau$, then ${a'_i}^\sigma = 
{a'_j}^\tau$, and, as $a_i'$ is contained in the tuple $\bar a$, also  
${a'_i}^\sigma = {a'_i}^\tau$. Because $\tau$ is an embedding, hence 
injective, we get $i=j$. Now $\bar a^\sigma = \bar a^\tau$ implies
$b_i^{f(\sigma)} = b_i^{f(\tau)}$ by the construction of $\bar a$, proofing 
$\rck{f}$ to be well defined. The surjectivity is clear by the choice of the 
$\bar b_i$. 
\ENDCLAIM

\medskip
\textsc{Claim:}
$\rck{f}$ is an existential interpretation.

Let $X \subseteq B^l$ and consider $\rck{f}^{-1}[X] \subseteq A^{(m+1)l}$. 
An element $\bar y$ therein has the form 
$$\bar y \,=\, \big((a'_{i_1},\bar a)^{\sigma_1},\dots,(a'_{i_l},
  \bar a)^{\sigma_l}\big) \ \text{ with }\ \big( b_{i_1}^{f(\sigma_1)},
  \dots, b_{i_l}^{f(\sigma_l)} \big) \;\in\; X.$$
Let $\sigma\in\Emb(\mf A)$. Then 
\begin{align*} 
  \bar y^\sigma \,=&\; \big( (a'_{i_1}, \bar a)^{\sigma_1\sigma}, \dots, 
  (a'_{i_l},\bar a)^{\sigma_l\sigma} \big) \;\in\; \rck{f}^{-1}[X] \\
  &\iff
  \big( b_{i_1}^{f(\sigma_1\sigma)}, \dots, b_{i_l}^{f(\sigma_l\sigma)} \big)  
  \,=\, \big( b_{i_1}^{f(\sigma_1)}, \dots, b_{i_l}^{f(\sigma_l)} 
  \big)_{}^{f(\sigma)} \;\in\, X.
\end{align*}
If $X$ is a basic set of the structure $\mf B$, then $X$ is closed under 
$\End(\mf B)$, thus the second condition is satisfied, whence 
$\rck{f}^{-1}[X]$ is existentially definable by Proposition~\ref{S-BP}.
\ENDCLAIM

``$\Longrightarrow$'':
Let $\mf B$ be existentially interpreted in $\mf A$ by a surjection 
$i: A^l \supseteq U \surj B$. Then $U = i^{-1}[B]$ and $E := i^{-1}[=_B]$ are 
existentially definable in $\mf A$, hence closed under $\Emb(\mf A)$ by  
Proposition~\ref{S-BP}. It follows that every $\sigma \in \Emb(\mf A)$ induces 
a map $\sigma^*: B \to B, uE \mapsto u^\sigma E$. Because the inverse image 
$i^{-1}[R]$ of every basic set $R$ of $\mf B$ is also existentially definable 
and thus closed under $\Emb(\mf A)$, the map $\sigma^*$ is even a 
homomorphism. This defines a mapping $\Emb(i): \Emb(\mf A) \to \End(\mf B),
\sigma \mapsto \sigma^*$, which clearly is a monoid homomorphism. The 
homomorphism is continuous: if $\bar a,\bar a'$ are inverse images of 
$\bar b, \bar b'$, then the open set $U_{\bar a,\bar a'}$ lies in the inverse 
image of $U_{\bar b,\bar b'}$. By the general theory developed in \cite{AZ}, 
the interpretation $i$ induces a continuous group homomorphism $\Aut(i): 
\Aut(\mf A) \to \Aut(\mf B)$ in the same way as above, that is $\Aut(i)$ is 
the map induced by $\Emb(i)$. Then by Theorem~1.2 in \cite{AZ}, $B$ is covered 
by finitely many orbits under the image of $\Aut(\mf A)$ under $\Aut(i)$.
\END

\begin{Rem}
In ``$\Longleftarrow$'', it follows in particular that $\mf B$ is  
$\aleph_0$-categorical, too. 

If this is known before, and if $\Aut(\mf B)$ is contained in the image of $f$ 
---in particular if $f$ is surjective--- then $B$ is automatically covered by 
finitely many orbits under the image of $f$. 

If the image of $f$ is contained in $\Emb(\mf B)$, then Proposition~\ref{S-BP} 
can be applied to $\mf B$, and the same argument as above shows that not only 
the inverse images of the basic sets, but of all existentially defined sets of 
$\mf B$ are existentially definable in $\mf A$.
\end{Rem}

\subsection{The positive existential case} \label{S-posex}

The proof of Theorem~\ref{Th-exist} works as well if one replaces ``$\Emb$'' 
by ``$\End$'' and ``existential'' by ``positive existential'', except for the 
well definedness of $\rck{f}$. The remark after Lemma~\ref{L-contractible} 
will show that there is no general solution to this problem. Therefore, we 
have to restrict our attention to a well behaved class of structures.

\begin{Def}
  An $\aleph_0$-categorical structure is called \emph{contractible} if it has 
  a constant endomorphism. 
\end{Def}

\begin{Lemma} \label{L-zus}
  An $\aleph_0$-categorical structure $\mf A$ is contractible if and only if 
  for each two tuples $\bar c_0,\bar c_1$ out of $A$ of same length there is 
  an endomorphism $\sigma \in \End(\mf A)$ such that $\bar c_0^\sigma = 
  \bar c_1^\sigma$. 
\end{Lemma} 

\PROOF 
Clearly, a contractible structure satisfies the condition. Assume now that the 
condition is satisfied. Given a tuple $\bar a = (a_1,\dots,a_k)$, choose an 
endomorphism $\sigma$ with $(a_1,\dots,a_k)^\sigma = (a_1,\dots,a_1)^\sigma$. 
Then $\sigma$ is constant $=c$ on $\bar a$. By multiplying with automorphisms 
we can assume that $c$ is an element of a fixed representation system 
$\{c_1,\dots,c_l\}$ of the $1$-types and moreover that it only depends on 
the type of $\bar a$. If we do this for a long tuple composed from 
representations of all $k$-types, we see that we can choose the value for each 
$k$ even independently from the type of $\bar a$. But then one of the finitely 
many values in question $c_1,\dots,c_l$ must work for every finite tuple, say 
$c$. Now $\End(\mf A)$ is closed ${}^A\!A$, therefore the constant map $c$ is 
an endomorphism. 
\END

\begin{Thm} \label{Th-posexist} 
  Let $\mf A$ be an $\aleph_0$-categorical, non-contractible structure. Then 
  $\mf B$ is positive existentially interpretable in $\mf A$ if and only if 
  there is a continuous monoid homomorphism $f:\End(\mf A)\to\End(\mf B)$
  such that $B$ is covered by finitely many orbits under the image of $f$, or, 
  equivalently, such that $B$ is covered by finitely many orbits under 
  $f[\Aut(\mf A)]$.
\end{Thm}

\PROOF
Take the proof of Theorem~\ref{Th-exist}, replace ``$\Emb$'' by
``$\End$'' and ``existential'' by ``positive existential'', and change
the definition of $U$ as follows: Choose tuples $\bar c_0,\bar c_1$ of
length $l$ such that $\bar c_0^\sigma \neq \bar c_1^\sigma$ for all
endomorphisms $\sigma$ as given by Lemma~\ref{L-zus}. Then let $a_i'$
be the $ml$-tuple  
$(\bar c_0,\dots,\bar c_0,\bar c_1,\bar c_0,\dots,\bar c_0)$ where
$\bar c_1$ is at the $i$th position. We may assume $m \geqslant 3$.  

Now if $(a'_i,\bar a)^\sigma = (a'_j,\bar a)^\tau$ for $i \neq j$,
then by comparing the appropriate coordinates we get $\bar c_1^\sigma =
\bar c_0^\tau = \bar c_0^\sigma$: contradiction. Thus again $\rck{f}$ is well
defined and everything goes through as in the proof of Theorem~\ref{Th-exist}. 
\END

\medskip
In fact, for the direction ``$\Longrightarrow$'' we do not need $\mf A$ to be 
non-contractible. Therefore:

\begin{Satz} \label{S-tsixesop}
  If $\mf A$ is an $\aleph_0$-categorical structure and if $i$ is a positive 
  existential interpretation of $\mf B$ in $\mf A$, then there is a continuous 
  monoid homomorphism $\End(i): \End(\mf A) \to \End(\mf B)$ such that $B$ 
  is covered by finitely many orbits under $f[\Aut(\mf A)]$. 
\end{Satz}

\begin{Cor}
``$\End$'' is a functor from the category of $\aleph_0$-categorical structures 
together with positive existential interpretations as morphisms into the 
category of topological monoids with continuous monoid homomorphisms. 
\end{Cor} 

\PROOF 
Check that the composition of positive existential interpretations is again 
positive existential (replacing a quantifier-free sub-formula of a positive 
existential formula by a positive existential formula yields again a positive 
existential formula). The rest follows from Proposition~\ref{S-tsixesop}.
\END

\medskip
Finally we remark that Theorem~\ref{Th-posexist} can't be extended to 
arbitrary $\aleph_0$-categorical structures:

\begin{Lemma} \label{L-contractible}
  If $\mf A$ is contractible and if $\mf B$ is positive existentially 
  interpretable in $\mf A$, then $\mf B$ is contractible, too.
\end{Lemma} 

\PROOF 
Let $\sigma$ be a constant endomorphism of $\mf A$ and let $\mf B$ be 
positive existentially interpreted in $\mf A$ by the interpretation $i$. Then 
$\sigma^* = \End(i)(\sigma)$ is a constant endomorphism of $\mf B$. 
\END

Note that there are non-contractible finite structures $\mf B$, which by 
Lemma~\ref{L-contractible} are not positive existentially interpretable in a 
contractible structure as for example $(\mathbb N,=)$, but the conditions of 
Theorem~\ref{Th-posexist} are trivially satisfied: the trivial monoid
homomorphism $\End(\mathbb N,=) \to \End(\mf B)$ is continuous and
$B$, being finite, is covered by finitely many orbits of the image
$\{\id\}$.

\subsection{Bi-interpretability}

For non-contractible structures, the theory of bi-interpretability of 
\cite{AZ} can be extended to positive existential interpretations. In the 
general case, or for existential interpretations, only partial results hold.

\begin{Def} \label{Def-homotop}
Following \cite{AZ}, we call two interpretations $i_1$ and $i_2$ of
$\mf B$ in $\mf A$ \emph{$\End$-homotopic} if $\End(i_1) = \End(i_2)$. 
Two $\aleph_0$-categorical structures $\mf A$ and $\mf B$ are
\emph{positive existentially bi-interpretable} if there are mutual 
positive existential interpretations $i$ and $j$ such that $i \circ j$
and $j \circ i$ are $\End$-homotopic to the identical interpretations
$\id_{\mf A}, \id_{\mf B}$ respectively.
\end{Def}

\begin{Lemma} \label{L-bi} 
  Let $i_1,i_2$ be two interpretations of $\mf B$ in $\mf A$. Then $\End(i_1) 
  = \End(i_2)$ holds if and only if the set $I_{i_1,i_2}:=\{(\bar x, \bar y) 
  \mid i_1(\bar x) = i_2(\bar y) \}$ is positive existentially definable in 
  $\mf A$. 
\end{Lemma} 
  
\PROOF
$\End(i_1), End(i_2)$ associate with an endomorphism $\sigma \in \End(\mf A)$ 
the maps induced by $\bar x\mapsto\bar x^\sigma, \bar y\mapsto\bar y^\sigma$, 
respectively. Both are the same if and only if $(\bar x,\bar y)\in I_{i_1,i_2}$ 
implies $(\bar x,\bar y)^\sigma = (\bar x^\sigma, \bar y^\sigma) \in 
I_{i_1,i_2}$. But according to Proposition~\ref{S-BP} this is exactly the case 
if $I_{i_1,i_2}$ is positive existentially definable. 
\END

\begin{Satz} \label{Pr-posexbi} 
  Let $\mf A$ and $\mf B$ be $\aleph_0$-categorical structures. If they are 
  positive existentially bi-interpretable, then $\End(\mf A)$ and $\End(\mf B)$ 
  are isomorphic as topological monoids. The converse holds for 
  non-contractible structures.
\end{Satz} 

\PROOF 
``$\Longrightarrow$'':
Let $i$ and $j$ be mutual interpretations witnessing the positive existential 
bi-interpretability. Then $j \circ i$ is $\End$-homotopic to the identical 
interpretation, hence $\End(j)\circ\End(i) = \End(j\circ i) = \End(\id) =\id$. 
Symmetrically, $\End(i)\circ\End(j) = \id$, hence $\End(i) = \End(j)^{-1}$ has 
to be a bi-continuous isomorphism. 

``$\Longleftarrow$'':
If $f: \End(\mf A) \to \End(\mf B)$ is an isomorphism, then by 
Theorem~\ref{Th-posexist}, $f$ and $f^{-1}$ yield interpretations $\rck f$ 
and $\rck f^{-1}$. The composition $j := \rck f^{-1}\circ\rck f:\mf A\to\mf A$ 
then induces the map $\End(\mf A)\to\End(\mf B)\to\End(\mf A), \sigma\mapsto 
f^{-1}(\sigma) \mapsto f(f^{-1}(\sigma)) = \sigma$, hence $\End(j) = \id = 
\End(\id)$. By symmetry, also $\End(\rck f\circ\rck f^{-1}) = \End(\id)$. 
\END

\medskip
The converse of Proposition~\ref{Pr-posexbi} does not hold for arbitrary 
$\aleph_0$-categorical structures, as Lemma~\ref{L-contractible} together with 
the following lemma shows.

\begin{Lemma} \label{L-isocontr}
  The isomorphism type of $\End(\mf A)$ does not determine whether
  $\mf A$ is contractible. 
\end{Lemma}

\PROOF
Let $\mf A$ be contractible $\mc L$-structure. We may assume the
language $\mc L$ to be relational. Let $\mf B$ be an $\mc L \cup
\{c,P\}$-structure that results from joining a new element $c$ to 
$\mf A$ and a predicate $P$ for the set $A$. Then $\mf B$ is not
contractible, but clearly $\End(\mf A)$ and $\End(\mf B)$ are
isomorphic. 
\END

\medskip 
On the other hand, each contractible structure contains an absorbing
endomorphism $\sigma$, i.e. $\tau\sigma = \sigma$ for every $\tau$ (and if 
there are constant endomorphisms, then they are exactly the absorbing 
elements). So non-contractibility can sometimes be seen from the endomorphism
monoid. 

Whether there are similar interpretability results for contractible structures 
is unclear.

\medskip
We will see in Section~\ref{S-Ex} that ``$\Emb$'' is not a functor as 
``$\Aut$'' and ``$\End$'' are. Therefore, the characterisation of existential 
bi-interpretability via the embedding monoids only holds in one direction.
With definitions analogously to Definition~\ref{Def-homotop} and the same 
proofs as for Lemma~\ref{L-bi} and Proposition~\ref{Pr-posexbi}, we get:

\begin{Satz} 
  Let $i_1,i_2$ be two interpretations of $\mf B$ in $\mf A$. Then $\Emb(i_1) 
  = \Emb(i_2)$ holds if and only if the set $I_{i_1,i_2}:=\{(\bar x, \bar y) 
  \mid i_1(\bar x) = i_2(\bar y) \}$ is existentially definable in $\mf A$. 
  
  If $\Emb(\mf A)$ and $\Emb(\mf B)$ are isomorphic as topological monoids, 
  then $\mf A$ and $\mf B$ are existentially  bi-interpretable.
\end{Satz}

The converse of the second part does not hold, as Example~\ref{Bsp5} below 
shows.

\section{Examples} \label{S-Ex}

We have seen that ``$\End$'' can be considered as a functor of the category of 
$\aleph_0$-categorical structures with positive existential interpretations 
into the category of topological monoids with continuous monoid homomorphisms.
This is not possible for ``$\Emb$'' and existential interpretations, for at 
least two reasons: $\aleph_0$-categorical structures with existential 
interpretations do not form a category, and the natural way to define $\Emb$ 
on morphisms leads to non-embeddings. We start with an example for the second 
problem:

\begin{Bsp} \label{Bsp2}
The image of a monoid homomorphism $f: \Emb(\mf A) \to \End(\mf B)$ is not in 
general contained in $\Emb(\mf B)$.

Let $\mf M_1$ be the following structure: $M_1$ is a countably infinite set, 
$E$ an equivalence relation on $M_1$ with infinitely many one-element classes, 
infinitely many two-element classes and no others. The language just contains 
a symbol for $E$. In $\mf M_1$, the structure $\mf M_2$ of an infinite, 
co-infinite predicate $P$ is existentially definable as $M_1/E$ with $P$ being 
the image of the two-element classes. Now there are embeddings of $\mf M_1$ 
mapping one-element classes into two-element classes. Their image in $\mf M_2$ 
are endomorphisms that are not strong.\END
\end{Bsp}

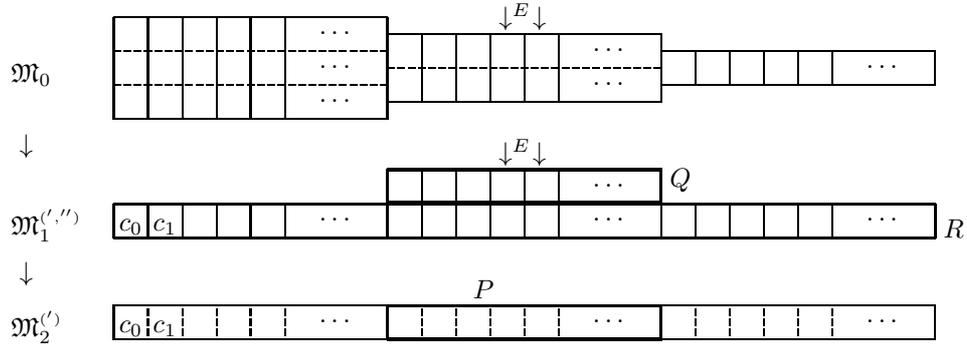
\begin{figure}[ht] \label{Abb}
\begin{center}
\unitlength.45cm
\begin{picture}(28,10)
\put(11,4){\line(1,0){8}}
\multiput(3,7.5)(0.3,0){27}{\line(1,0){0.2}}
\multiput(3,8.5)(0.3,0){27}{\line(1,0){0.2}}
\multiput(11,8)(0.3,0){27}{\line(1,0){0.2}}

\put(0,0.1){$\mf M_2^{\scriptscriptstyle(')}$}
\put(3,0){\line(1,0){8}}
\put(3,1){\line(1,0){8}}
\put(19,0){\line(1,0){8}}
\put(19,1){\line(1,0){8}}
\put(3,0){\line(0,1){1}}
\put(27,0){\line(0,1){1}}
\put(13.5,1.2){$P$}
\put(3.15,0.2){$c_0$}
\put(4.15,0.2){$c_1$}
\put(9,.5){$\ldots$}
\put(17,.5){$\ldots$}
\put(25,.5){$\ldots$}
\multiput(3,0.1)(1,0){6}{\multiput(0,0)(0,0.3){3}{\line(0,1){0.2}}}
\multiput(11,0.1)(1,0){6}{\multiput(0,0)(0,0.3){3}{\line(0,1){0.2}}}
\multiput(19,0.1)(1,0){6}{\multiput(0,0)(0,0.3){3}{\line(0,1){0.2}}}

\thicklines
\put(11,0){\line(1,0){8}}
\put(11,1){\line(1,0){8}}
\put(11,0){\line(0,1){1}}
\put(19,0){\line(0,1){1}}

\thinlines 
\put(0.25,1.75){$\downarrow$}
\put(0,3.1){$\mf M_1^{\scriptscriptstyle(','')}$}
\put(9,3.5){$\ldots$}
\put(17,3.5){$\ldots$}
\put(17,4.5){$\ldots$}
\put(25,3.5){$\ldots$}
\put(27.25,3){$R$}
\put(19.25,4.5){$Q$}
\put(3.15,3.2){$c_0$}
\put(4.15,3.2){$c_1$}
\multiput(3,3)(1,0){6}{\line(0,1){1}}
\put(11,3){\line(0,1){1}}
\multiput(11,3)(1,0){6}{\line(0,1){2}}
\put(19,3){\line(0,1){2}}
\multiput(19,3)(1,0){6}{\line(0,1){1}}
\put(27,3){\line(0,1){1}}
\put(14.25,5.25){$\downarrow^E$}
\put(15.25,5.25){$\downarrow$}
\thicklines
\put(3,3){\line(1,0){24}}
\put(3,4){\line(1,0){24}}
\put(11,4.05){\line(1,0){8}}
\put(11,5){\line(1,0){8}}
\put(3,3){\line(0,1){1}}
\put(11,4){\line(0,1){1}}
\put(19,4){\line(0,1){1}}
\put(27,3){\line(0,1){1}}

\put(0.25,5.5){$\downarrow$}
\put(0,7.6){$\mf M_0$}
\thinlines
\put(3,6.5){\line(1,0){8}}
\put(3,9.5){\line(1,0){8}}
\put(11,7){\line(1,0){8}}
\put(11,9){\line(1,0){8}}
\put(19,7.5){\line(1,0){8}}
\put(19,8.5){\line(1,0){8}}
\put(3,6.5){\line(0,1){3}}
\put(11,6.5){\line(0,1){3}}
\put(19,7){\line(0,1){2}}
\put(27,7.5){\line(0,1){1}}
\multiput(3,6.5)(1,0){6}{\line(0,1){3}}
\multiput(11,7)(1,0){6}{\line(0,1){2}}
\multiput(19,7.5)(1,0){6}{\line(0,1){1}}
\put(9,7){$\ldots$}
\put(9,8){$\ldots$}
\put(9,9){$\ldots$}
\put(17,7.5){$\ldots$}
\put(17,8.5){$\ldots$}
\put(25,8){$\ldots$}
\put(14.25,9.25){$\downarrow^E$}
\put(15.25,9.25){$\downarrow$}
\end{picture}

\medskip
Squares correspond to elements of the structures; dotted lines do not 
correspond to structure named in the signature.

\caption{Examples~\ref{Bsp2}, \ref{Bsp4}, \ref{Bsp3} and \ref{Bsp5}.}
\end{center}
\end{figure}

This phenomenon is in connection with the following: $\Aut(\mf M)$ can
be characterised in the abstract monoid $\End(\mf M)$ as the subgroup
of invertible elements $\End(\mf M)^*$. Therefore a homomorphism
between endomorphism monoids restricts to a homomorphism between the
automorphism group. $\Emb(\mf M)$, on the other hand, can only be
defined in the ``permutation monoid''\footnote{I.e. the monoid with
  its action on the set $M$; in analogy to ``permutation group'', though 
  the elements of the monoid are not in generally acting as permutations.} 
$\End(\mf M)$; no characterisation
in the abstract monoid is possible as Example~\ref{Bsp3} shows.  

The following holds in general, with $E = \End(\mf M)$:
$$\{\sigma \in E \mid \exists \tau \in E\; \sigma\tau \in E^* \}
  \subseteq \Emb(\mf M) \subseteq \{\sigma \in E \mid \forall
  \tau_1,\tau_2 \in E\; (\tau_1\sigma = \tau_2\sigma \Rightarrow \tau_1 = \tau_2)\}$$

\begin{Bsp} \label{Bsp4}
Take $\mf M_1,\mf M_2$ as in Example~\ref{Bsp2} and expand $\mf M_1$ to the 
structure $\mf M_1''$ by adding a predicate $Q$ that picks exactly one 
element out of each two-element class, and a predicate $R$ for its complement.
It is easy to see that the interpretation of $\mf M_2$ in $\mf M_1''$
induces an isomorphism $\End(\mf M_2) \to \End(\mf M_1'')$. The image
of the injective endomorphisms of $\mf M_2$ are exactly the injective
endomorphisms of $\mf M_1''$, but the image of $\Emb(\mf M_2)$ is only
a proper subset of $\Emb(\mf M_1'')$. Thus ``$\Emb$'' does not allow
an abstract characterisation. 
(Note that $\mf M_2$ is contractible, but $\mf M_1$ is not, so they are not 
positive existentially bi-interpretable.)\END
\end{Bsp}

\begin{Bsp} \label{Bsp3}
The composition of existential interpretations need not to be existential:

Let $\mf M_0$ be the structure on a countably infinite set $M_0$ with an 
equivalence relation with infinitely many three-element, two-element and 
one-element classes and no others. It interprets existentially the structure 
$\mf M_1$ in Example~\ref{Bsp2} by collapsing each three-element class to one 
element. The composition with the interpretation in Example~\ref{Bsp2} however 
is not existential: it yields the interpretation of $\mf M_2$ in $\mf M_0$, 
which on $M_0/E$ define a predicate for the images of the two-element classes. 
This predicate is not existential.\END
\end{Bsp}

We conclude with an example of two structures $\mf M_1',\mf M_2'$ that are 
mutually positive existentially interpretable, existentially but not positive 
existentially bi-interpretable, and with non isomorphic embedding monoids.

\begin{Bsp} \label{Bsp5}
  $\mf M_1'$ is the structure $\mf M_1$ from Example~\ref{Bsp2} with
  an additional predicate $Q$, that takes exactly one element out of
  each two-element equivalence class, and with two additional
  non-equivalent constants $c_0,c_1$,living in one-element classes. 
  $\mf M_2'$ is the structure $\mf M_2$ of an infinite, co-infinite
  predicate $P$ from Example~\ref{Bsp2} together with two distinct
  constants $c_0,c_1$ not in $P$. 

  We interpret $\mf M_2'$ positive existentially in $\mf M_1'$ as
  $M_1/E$ with $\exists y (Qy \land Exy)$ providing the predicate $P$ and 
  by keeping the two constants. 
  We interpret $\mf M_1'$ positive existentially in $\mf M_2'$ as
  follows: the universe is 
  $(M_2 \times\{c_0\}) \cup (P \times\{c_1\})$;
  the equivalence relation $E$ is ``same first coordinate'', 
  the predicate $Q$ is $P \times\{c_1\}$ and the
  two constants are $(c_0,c_0)$ and $(c_1,c_0)$.

  Both structures are bi-interpretable as they have the same
  automorphism group $S_\omega \times S_\omega$. But the endomorphism
  monoids are not isomorphic: In $\End(\mf M_1')$, there is the
  endomorphism $\sigma$ that collapses all two-element classes and is
  the identity on $Q$ and the one-element classes. This endomorphism
  satisfies $\sigma^2 = \sigma$ and commutes with all automorphisms,
  hence is definable in the structure. There are three such elements
  in $\End(\mf M_2')$: identity on $P \cup \{c_0,c_1\}$ and either
  identity or constant $= c_0$ or $c_1$ on the rest, but
  six in $\End(\mf M_1')$: the corresponding maps and their
  compositions with $\sigma$. 

  It is easy to verify that the bi-interpretation above is in fact an 
  existential bi-interpretation. But $\Emb(\mf M_1')\not\isom \Emb(\mf M_2')$, 
  as can be seen with the following argument:  
  Because of the bi-interpretability, both structures have isomorphic 
  automorphism groups, of isomorphism type $S_\omega\times S_\omega$. If the 
  two embedding monoids were isomorphic, an isomorphism had to respect this 
  decomposition as it is unique in this group. 
  Now in $\mf M_2'$, each embedding $\sigma$ is a (commuting) product
  of the two embeddings $\sigma\!\!\restriction_{P} \cup\, \id_{M_2\ohne P}$
  and $\id_{P} \cup \sigma\!\!\restriction_{M_2\ohne P}$, and each of
  the two commutes with one of the factors $S_\omega$. In $\mf M_1$
  however, there are embeddings which move one-element equivalence
  classes into two-element classes. Such an endomorphism cannot be
  decomposed in that way.  
  \END
\end{Bsp}

\section{Concluding remarks} \label{S-Con}

We have shown characterisations of existential and positive existential 
interpretability in $\aleph_0$-categorical structures:
\begin{itemize}
\item 
  A structure $\mathfrak B$ has an existential interpretation in an 
  $\aleph_0$-categorical structure $\mathfrak A$  if and only if there is a 
  continuous monoid homomorphism $f$ from the monoid of self-embeddings of 
  $\mathfrak A$ to the endomorphism monoid of $\mathfrak B$ such that the 
  domain of $\mathfrak B$ is covered by finitely many orbits under the image 
  of $f$.
\item 
  A structure $\mathfrak B$ has a positive existential interpretation in an   
  $\aleph_0$-categorical structure $\mathfrak A$ without constant 
  endomorphisms if and only if there is a continuous monoid homomorphism $f$ 
  from the endomorphism monoid of $\mathfrak A$ to the endomorphism monoid of 
  $\mathfrak B$ such that the domain of $\mathfrak B$ is covered by finitely 
  many orbits under the image of $f$.
\end{itemize}

It is open whether the second result also holds for $\aleph_0$-categorical 
structures $\mathfrak{A,B}$ with constant endomorphisms.

It would be very interesting to find an analogous characterisation of 
\emph{primitive positive interpretability}. A formula is called 
\emph{primitive positive} if it is of the form 
$$\exists x_1 \dots \exists x_n (\psi_1 \wedge \dots \wedge \psi_m)$$ 
where $\psi_1,\dots,\psi_m$ are atomic formulas. Primitive positive 
interpretations play an important role for the study of the computational 
complexity of \emph{constraint satisfaction problems}. For a structure 
$\mathfrak A$ with finite relational signature $\tau$, the \emph{constraint 
satisfaction problem for $\mathfrak A$, CSP$(\mathfrak A)$,} is the 
computational problem to decide whether a given primitive positive 
$\tau$-sentence is true in $\mathfrak A$. Such problems are abundant in many 
areas of computer science.

It is well-known that if (every relation of) a structure $\mathfrak B$ is 
primitive positively definable in a structure $\mathfrak A$, then 
CSP$(\mathfrak B)$ has a polynomial-time reduction to CSP$(\mathfrak A)$. 
Indeed, an important technique to show that CSP$(\mathfrak A)$ is NP-hard is 
to find another structure $\mathfrak B$ such that CSP$(\mathfrak B)$ is 
already known to be NP-hard, and to give a primitive positive definition of 
$\mathfrak B$ in $\mathfrak A$. 

Primitive positive definability in an $\aleph_0$-categorical structure 
$\mathfrak A$ is captured by the \emph{polymorphisms} of $\mathfrak A$. 
A polymorphism of $\mf M$ is a homomorphism of some power $\mf M^n$ 
(with the product structure) to $\mf M$. A subset $X \subseteq M^k$ is called
\emph{closed under polymorphisms} if for all $n$, every polymorphism $\sigma: 
\mf M^n \to \mf M$ and all $\bar a_1,\dots,\bar a_n \in X$ we have 
$(\bar a_1,\dots,\bar a_n)^\sigma \in X$. The following has been shown
in \cite{BN}:

\begin{Thm} \label{R-Poly}
  Let $M$ be an $\aleph_0$-categorical structure and $X \subseteq M^k$. 
  Then $X$ is positive primitive modulo the theory of $\mf M$ if and only if 
  $X$ is closed under polymorphisms.  
\end{Thm} 

%
%

The classification of the computational complexity of CSP$(\mathfrak A)$ for 
all highly set-transitive structures $\mathfrak A$ obtained 
in~\cite{BodirskyK08} makes essential use of this theorem.

An even more powerful tool to classify the computational complexity of 
CSP$(\mathfrak A)$ is \emph{primitive positive interpretability}. It has been 
shown in \cite{Bodirsky08} that if a structure $\mathfrak B$ 
has a primitive positive interpretation in $\mathfrak A$, then there is a 
polynomial-time reduction from CSP$(\mathfrak B)$ to CSP$(\mathfrak A)$.
Hence, it would be interesting to have algebraic characterisations of 
primitive positive interpretability in $\aleph_0$-categorical structures.

Note that the set of all polymorphisms of a structure $\mathfrak A$ can be 
seen as an algebra whose operations are precisely the polymorphisms of 
$\mathfrak A$; we will refer to this algebra as the \emph{polymorphism 
clone} of $\mathfrak A$. In fact, the set of all polymorphisms forms an object 
called \emph{clone} in universal algebra. The following characterisation of 
primitive positive interpretability has also been given 
in~\cite{Bodirsky08}.

\begin{Thm}\label{thm:pp-interpret}
  Let $\mathfrak A$ be finite or $\aleph_0$-categorical. Then a structure 
  $\mathfrak B$ has a primitive positive interpretation in $\mathfrak A$ if 
  and only if  there is an algebra $\mathbb B$ in the pseudo-variety generated 
  by the polymorphism clone of $\mathfrak A$ such that all operations of 
  $\mathbb B$ are polymorphisms of $\mathfrak B$.
\end{Thm}

It follows that for $\aleph_0$-categorical structures $\mathfrak A$ the 
computational complexity of CSP$(\mathfrak A)$ is determined by the 
pseudo-variety generated by the polymorphism clone of $\mathfrak A$. 
We would like to give an alternative characterisation of primitive positive
interpretability in terms of the \emph{topological polymorphism clone} of 
$\mathfrak A$, in analogy to the theorems shown in this paper. In fact, we 
conjecture that the computational complexity of CSP$(\mathfrak A)$ is indeed 
determined by the topological polymorphism clone of $\mathfrak A$.

\bibliographystyle{plain}
\bibliography{endomorphisms} 

\end{document}